\documentclass{amsart}

\usepackage{graphicx}
\usepackage{epstopdf}
\usepackage{lineno}
\usepackage{amsrefs}
\usepackage[justification=raggedright]{caption}

\newtheorem{thm}{Theorem}[section]
\newtheorem{lem}[thm]{Lemma}
\newtheorem{prop}[thm]{Proposition}
\newtheorem{cor}[thm]{Corollary}
\newtheorem*{cor2}{Corollary}
\newtheorem*{thmA1}{Theorem A1}
\newtheorem*{thmA2}{Theorem A2}
\newtheorem*{thmB}{Theorem B}
\newtheorem*{thmC1}{Theorem C1}
\newtheorem*{thmC2}{Theorem C2}
{\theoremstyle{definition} 
\newtheorem{defn}[thm]{Definition}
\newtheorem{exmp}[thm]{Example}
\newtheorem{rmk}[thm]{Remark}}
\newtheorem{question}{Question}

\newcommand\bbn{\mathbb N}
\newcommand\bbz{\mathbb{Z}}
\newcommand\bbr{\mathbb{R}}
\newcommand\Om{\Omega}
\newcommand\Sig{\Sigma}
\newcommand\g{\gamma}
\newcommand\sig{\sigma}
\newcommand\M{\mathcal{M}}
\newcommand\B{\mathcal{B}}
\newcommand\A{\mathcal{A}}
  
\newcommand{\incf}[4]{
\begin{figure}[!ht]
\includegraphics[width=0.#1\textwidth]{./figures/#2}
\caption{#3}
\label{#4}
\end{figure}
}

\title[Construction of Knots as Murasugi Sums of Seifert Surfaces]{Construction and Decomposition of Knots as Murasugi Sums of Seifert Surfaces}

\author{Jared Able}
\address{Department of Mathematics, Michigan State University, East Lansing, Michigan 48824}
\email{ablejare@msu.edu}
\thanks{The first author was supported by NSF grant DMS-1709016}

\author{Mikami Hirasawa}
\address{Department of Mathematics, Nagoya Institute of Technology, Nagoya Aichi 466-8555 Japan}
\email{hirasawa.mikami@nitech.ac.jp}
\thanks{The second author was supported by Grant-in-Aid for Scientific Research (C)}


\begin{document}

\begin{abstract}
A fixed knot $K$ acts via Murasugi sum on the space $\mathcal{S}$ of isotopy classes of knots. 
This operation endows $\mathcal{S}$ with a directed graph structure denoted by $M\kern-1pt SG(K)$.
We show that any given family of knots in $M\kern-1pt SG(K)$ has
the structure of a bi-directed complete graph, which is not the case if we restrict 
the complexity of Murasugi sums.
For that purpose, we show that any knot is a Murasugi sum of
any two knots, and we give lower and upper bounds for the minimal complexity of Murasugi sum
to obtain $K_3$ by $K_1$ and $K_2$. As an application, we show that given any three knots,
there is a braid for one knot which splits along a string into braids for the other two knots.
\end{abstract}

\maketitle

\section{Introduction}\label{sec:intro}
One development in knot theory is
to define and study the structure of the topological space $\mathcal{S}$
composed of isotopy classes of knots.
Here, a knot is an embedding of an oriented circle in the 3-sphere $S^3$ which can be
represented as a finite number of line segments, and a link is a disjoint union of knots.
We identify knots which are isotopic to one another.

In \cite{hu}, the {\it Gordian Complex} $\mathcal{G}$ of knots was defined as follows:
The vertex set of $\mathcal{G}$ consists of isotopy classes of knots, and 
a set of $n+1$ vertices $K_{0}, \dots, K_{n}$ spans an $n$-simplex
if and only if any pair of knots in it can be changed into each other by a single crossing change.
Since then, many studies have been done by replacing
the crossing change with other local operations on knots and on virtual knots as in \cites{lyz,gmv,hkys,ik,no}. 

In the following, we consider the operation
of Murasugi sum along Seifert surfaces of links. An oriented, embedded surface $F$ without a closed component is 
called a
{\it Seifert surface} for an oriented link $L$ if its boundary $\partial F$ coincides with the oriented link $L$.
An $m$-Murasugi sum is an operation to glue two Seifert surfaces $F_1$ and $F_2$ along an $m$-gon with $m$ even
(for the precise definition, see Definition \ref{dfn:msum}). 
 
 We may regard a fixed knot $K$ as an operation on the space of knots as follows.
 A knot $K_1$ is changed to a knot $K_2$ via $K$ if a Seifert surface for $K_2$ is
 obtained from a Seifert surface $F_1$ for $K_1$ 
 by Murasugi-summing a Seifert surface $F$ for $K$ to $F_1$.
 Thus for each fixed $K$, the space of knots $\mathcal{S}$ has the structure of directed graph,
 where an edge is an arrow or a double-headed arrow.

 \begin{defn}\label{dfn:m-graph}
 For a knot $K$,
 the Murasugi sum graph of knots $M\kern-1pt SG(K)$ is 
 a directed graph such that
 (1) the vertex set consists of all isotopy classes of knots,
 (2) two vertices $K_1$ and $K_2$ are connected by an edge
 with arrowhead on $K_2$ if
 there exist Seifert surfaces $F, F_1, F_2$ for $K,K_1, K_2$ such
 that $F_2$ is a Murasugi sum of $F$ and $F_1$.
 The {\it restricted Murasugi sum graph} $M\kern-1ptSG(K, n)$ is 
 considered by only allowing Murasugi sums along $m$-gons
 with $m\leq n$.
 \end{defn}

For $n=2$, a $2$-Murasugi sum is the connected sum operation. 
Hence for any knot $K$, 
$M\kern-1ptSG(K, 2)$ has an obvious structure where the
edges are arrows from each knot $K'$ to the connected sum of 
 $K$ and $K'$. For $n=4$, a $4$-Murasugi sum is called a {\it plumbing}.
It was shown in \cites{mura1, gabai, stall} that nice geometric properties of knots and surfaces
(such as fiberedness and genus-minimality)
were preserved under Murasugi sums and decompositions
of minimal genus Seifert surfaces. On the other hand, Thompson \cite{thomp} gave examples where
the trefoil is obtained as a plumbing of two unknots, and the unknot is obtained as a plumbing of two figure-eight knots. Thus,
expectations to generalize preservation results to Murasugi sums of 
non-minimal Seifert surfaces were negated.

In this paper, we generalize Thompson's examples to show that given any three knots, we can produce one of them as a Murasugi sum of the other two.

\begin{thmA1}
For any three knots $K_1, K_2, K_3$, there exist
Seifert surfaces $F_1, F_2,$ $F_3$ for them such that
$F_3$ is a Murasugi sum of $F_1$ and $F_2$.
\end{thmA1}
Therefore we have the following:

\begin{cor2}
For a knot $K$, any set of knots
$\{K_1, K_2, \dots, K_n\}$ in
$M\kern-1ptSG(K)$
composes a complete graph where all the edges are bi-directed.
\end{cor2}

We refine the result of Theorem A1 by giving an algorithm to find a closed braid for $K_3$ which naturally splits into closed braids for $K_1$ and $K_2$. See Figures \ref{fig:ex1},\,\ref{fig:ex2} following Example \ref{ex:knots-braid}.
\begin{thmA2}
For any three knots $K_1, K_2, K_3$, there are braids $b_1, b_2, b_3$ such that
$K_i$ is the closure of $b_i$  $(i=1,2,3)$, satisfying
the following:

If the braid $b_3$ is expressed as a braid word $W_3$
with generators\\
$\sigma_1, \sigma_2, \dots, \sigma_k, \dots, \sigma_n$,
then $W_1$ (resp. $W_2$) is obtained from
$W_3$ by deleting the generators $\sigma_1, \dots, \sigma_k$ (resp. $\sigma_{k+1}, \dots, \sigma_n$).
\end{thmA2}

To further study the structure of $M\kern-1ptSG(K, n)$ where
the size of Murasugi sums is limited, we give, in Section
\ref{sec:ngon}, lower and upper bounds on the minimal $m$-gon required to form $K_3$ as a Murasugi sum of $K_1$ and $K_2$. Our bounds are in terms of $d_{cb}(K,K')$ and $d_{bt}(K,K')$, which are, respectively, the minimal number of coherent band surgeries (resp. band-twists) required to transform $K$ into $K'$. For the precise definition of these operations, see Definitions \ref{dfn:bandtwist} and \ref{dfn:bandsurg}.

\begin{thmB}
If $K_3$ is an $m$-Murasugi sum of $K_1$ and $K_2$ such that $m$ is minimal, then 
\begin{center}
$d_{cb}(K_1\# K_2, K_3)+2\leq m\leq 2(d_{bt}(K_1,K_3)+d_{bt}(K_2,O)+1), $
\end{center}
where the roles of $K_1$ and $K_2$ can be switched to improve the upper bound. 
\end{thmB}

Another motivation for the study of manipulation of Seifert surfaces
is to provide tools for physically constructing knots as in
DNA, polymers, and proteins.
Manipulating Seifert surfaces seems handier than manipulating the knot itself,
while dissolving the surface yields a knot on the boundary. So in the Appendix (Section \ref{sec:append}), we 
prove the following two results regarding particularly nice Murasugi decompositions of Seifert surfaces of links in terms of positive Hopf bands, planar surfaces, and punctured Heegaard surfaces (i.e., the boundary of a standardly embedded handlebody).

\begin{thmC1}\label{thm:heggaard-positive}
Any link $L$ has a Seifert surface $F$ which is
a Murasugi sum of 
a once-punctured Heegaard surface
and
a boundary connected sum of positive Hopf bands.
\end{thmC1}

\begin{thmC2}\label{thm:three-planar}
For any link $L$, there exists a Seifert surface $F$ for
$L$ which consists of three planar surfaces 
$P, P_1, P_2$, where $P_1$ (resp. $P_2$) is Murasugi summed
to $P$ on the positive (resp. negative) side.
\end{thmC2}

\section{Any knot is a Murasugi sum of any two knots}\label{sec:anyknotsum}

The original construction of the Murasugi sum was first introduced by Murasugi in \cite{mura1} and was later coined as the Murasugi sum by Gabai in \cite{gabai}. For simplicity, we define the Murasugi sum in terms of Murasugi decomposition. See Figure \ref{fig:6sum}.
\begin{defn}\label{dfn:msum}
Let $F$ be a Seifert surface in $S^3$, let $\Sig$ be a 2-sphere such that $S^3\setminus \Sig$ is a union of two open 3-balls $B_1,B_2$ and such that $\Sig \cap F$ is an $m$-gon $\Om$ with $m$ even. Denote $F\cap \overline{B_1}=F_1$ and $F\cap \overline{B_2}=F_2$. Then we say that $F$ {\it decomposes} into $F_1$ and $F_2$ along $\Om$.
\end{defn}

\incf{50}{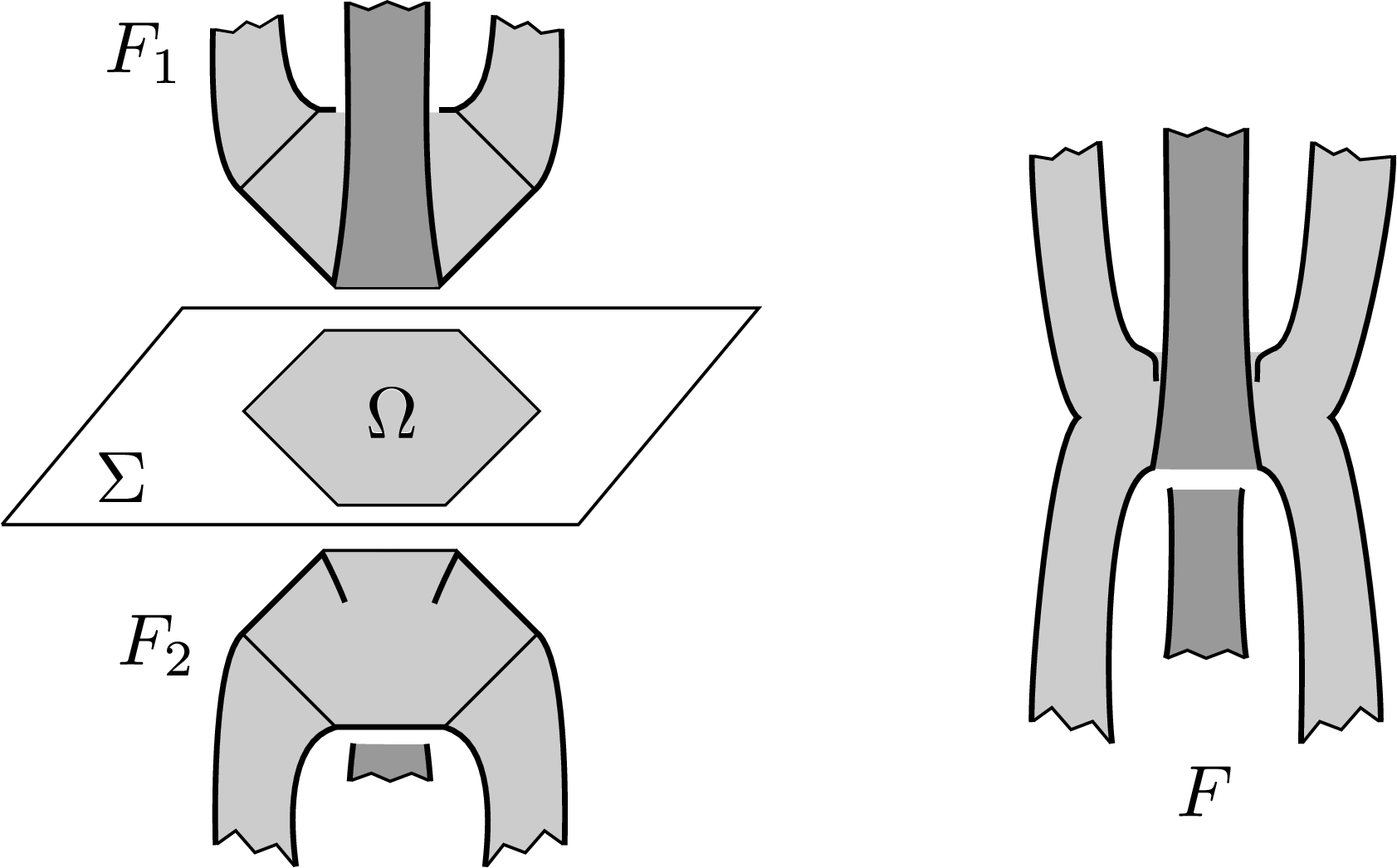}{A local picture of a 6-Murasugi sum}{fig:6sum}

If $F$ decomposes into $F_1$ and $F_2$ along an $m$-gon, then $F$ is said to be an $m$-Murasugi sum of $F_1$ and $F_2$, which we denote by $F_1 \star_m F_2$.
Given two knots $K_1,K_2$, we write $K_1\star_m K_2$ to denote the boundary of $F_1 \star_m F_2$ for some Seifert surfaces $F_1,F_2$ for $K_1,K_2$. Note that the summing disk $\Om$ can initially appear stretched and twisted, but we can isotope the given surfaces to see $\Om$ as flat. We can also isotope the surfaces so that $\Sig$ corresponds to the plane $z=0$ in $\bbr^3$, and in this situation, if $F_1$ lies above (resp. below) $z=0$ and $F_2$ lies below (resp. above) $z=0$, then we say that we Murasugi sum $F_1$ onto the positive (resp. negative) side of $F_2$.

There are several operations one can perform on knot diagrams to obtain a new knot. One such operation is a crossing change, and more generally, an antiparallel full-twisting.

\begin{defn}\label{dfn:bandtwist}
An {\it antiparallel full-twisting} on an oriented link is a local move where we
select a pair of
locally antiparallel strings 
 \raisebox{-5pt}{\includegraphics[height=13pt]{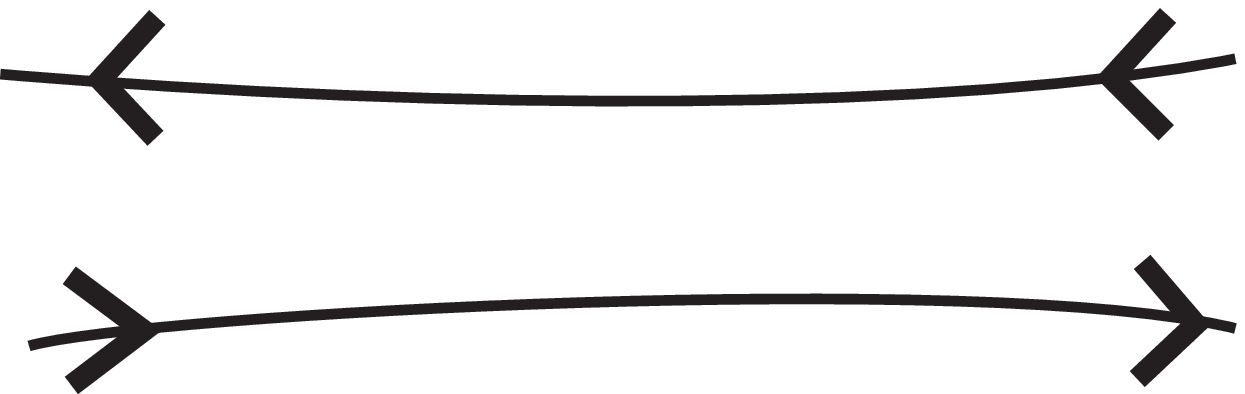}} and apply some number of full twists
  \raisebox{-5pt}{\includegraphics[height=13pt]{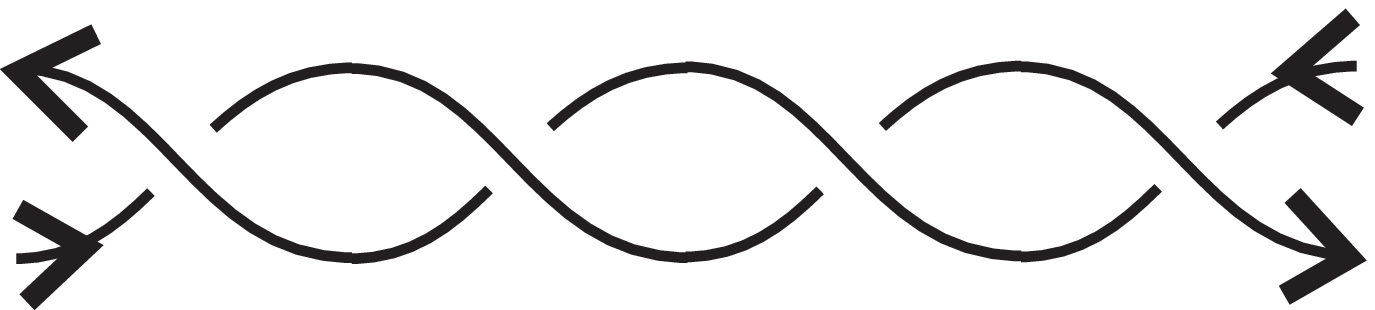}}.
 \end{defn}

For convenience, we sometimes refer to antiparallel full-twisting as band twisting. We can realize this twisting operation along an arc $\alpha$, 
where $\alpha$ is a short, unknotted arc connecting two antiparallel strings of a link $L$, which is contained within 
a small ball $B$ such that $L\cap B$ is a trivial $2$-string tangle.
In this setting, we can span a Seifert surface $F$ for $L$ such that
 $F \cap B$ is a rectangle $b$ containing the arc $\alpha$.
Then the twisting operation is realized by applying some full-twists to $b$.

Consider the two Seifert surfaces of the unknot in Figure \ref{fig:twounknots}. The following proposition states that a crossing change within a knot can be realized by either plumbing or deplumbing these surfaces. More generally, by increasing the number of full twists in $R_+$ or $R_-$, we can realize any antiparallel full-twisting by either plumbing or deplumbing Seifert surfaces for the unknot. 

 \incf{32}{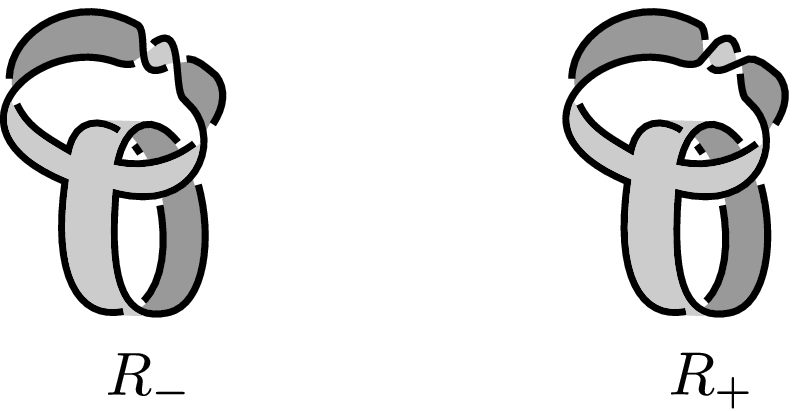}{Two Seifert surfaces for the unknot}{fig:twounknots}

\begin{prop}\label{prop:cc-by-plumbing}
Let $K_1$, $K_2$ be knots such that
$K_2$ is obtained by changing a positive crossing in $K_1$.
Then there exist Seifert surfaces $F_1$, $F_1'$ for $K_1$,
$F_2$, $F_2'$ for $K_2$, and $R_{+},R_{-}$ for the unknot
satisfying the following:
\begin{enumerate}
\item $F_1'$ is a plumbing of $F_2$ and $R_{+}$.
\item $F_2'$ is a plumbing of $F_1$ and $R_{-}$.
\end{enumerate}
\end{prop}

\begin{proof}
We illustrate both statements in Figure \ref{fig:crosschange}.
\incf{90}{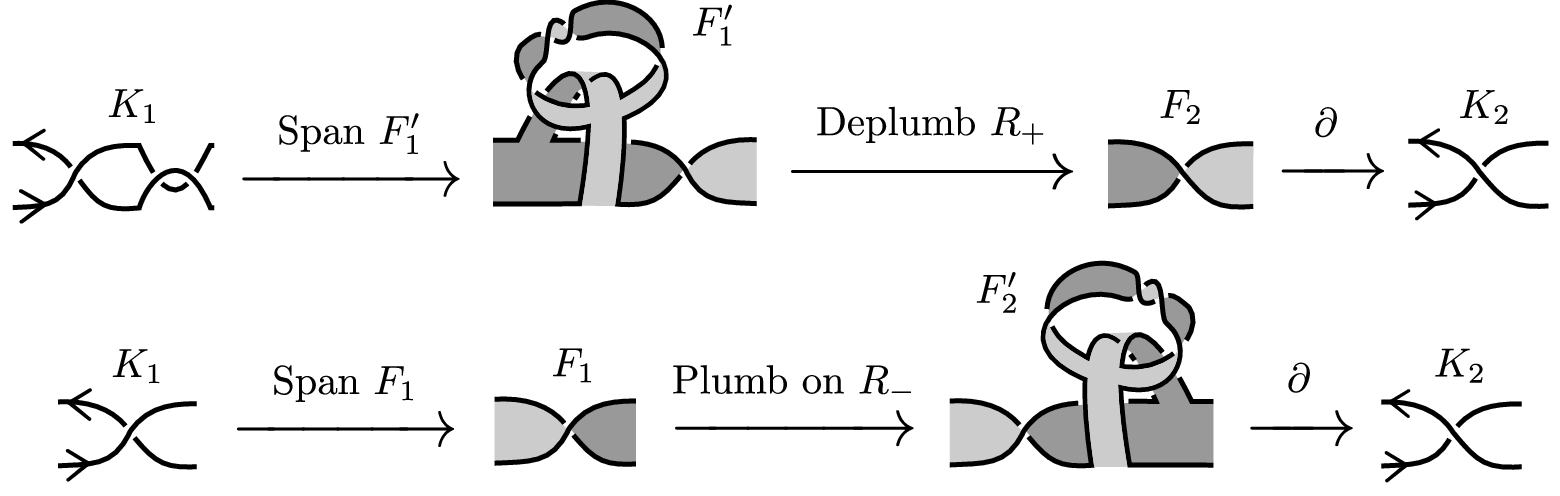}{Changing a positive crossing}{fig:crosschange}
\end{proof}
Furthermore, we can perform any number of crossing changes simultaneously via a single Murasugi sum with a Seifert surface for the unknot, 
which allows us to prove the following lemma.

 \begin{lem}\label{lemm:sum-of-unknots}
 Any knot $K$ has a Seifert surface $F$ 
 which is a Murasugi sum of 
 two Seifert surfaces $F_1, F_2$ for the unknot.
 \end{lem}
 
\begin{proof}
 For any diagram $D$ of $K$, we can choose a
 subset ${\mathcal C}$ of crossings such that
 we obtain the unknot by simultaneously changing
 the crossings in ${\mathcal C}$.
 To see this, start at some point in $D$ and walk along
 the knot. Then can specify ${\mathcal C}$ to consist of  
 those crossings which we enter first along an under-path and later along an over-path.
 Near each crossing
\raisebox{-5pt}{\includegraphics[height=13pt]{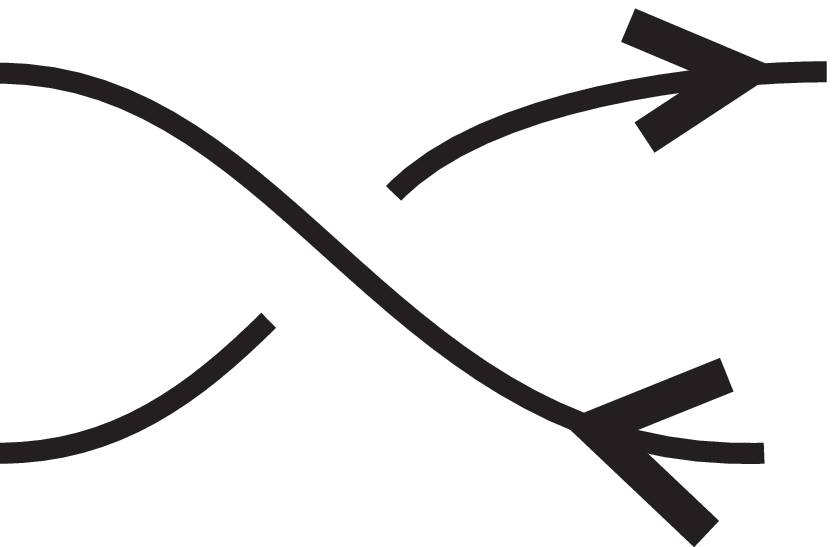}}
  in ${\mathcal C}$, 
 apply a Reidemeister II move to introduce an
 antiparallel clasp
 \raisebox{-5pt}{\includegraphics[height=13pt]{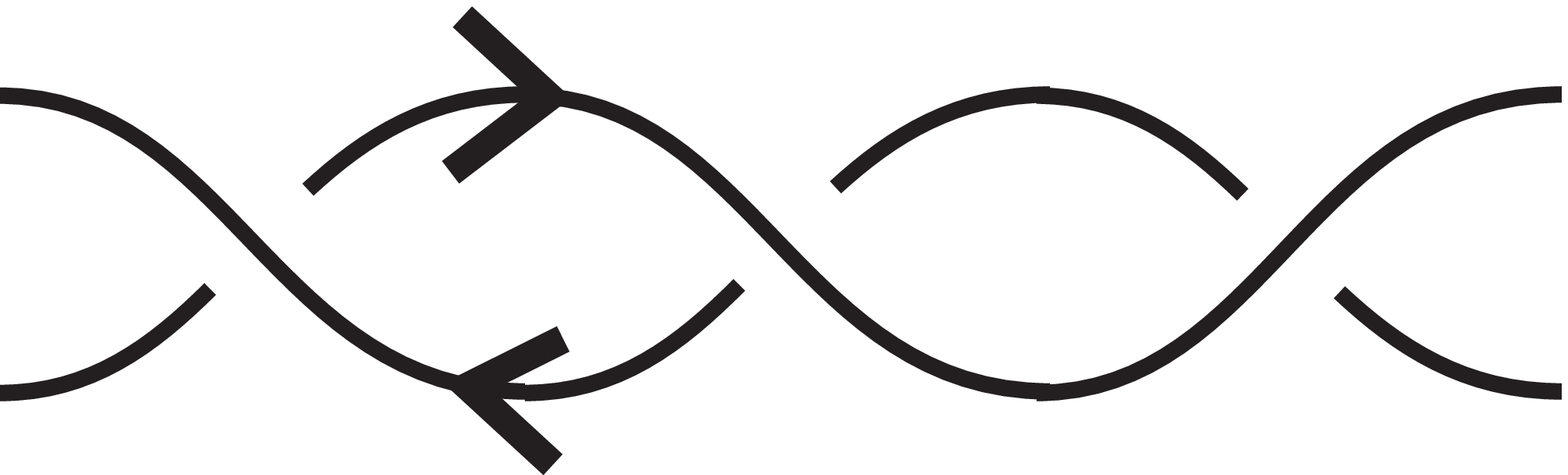}}, where undoing the clasp results in changing that crossing as in  \raisebox{-5pt}{\includegraphics[height=13pt]{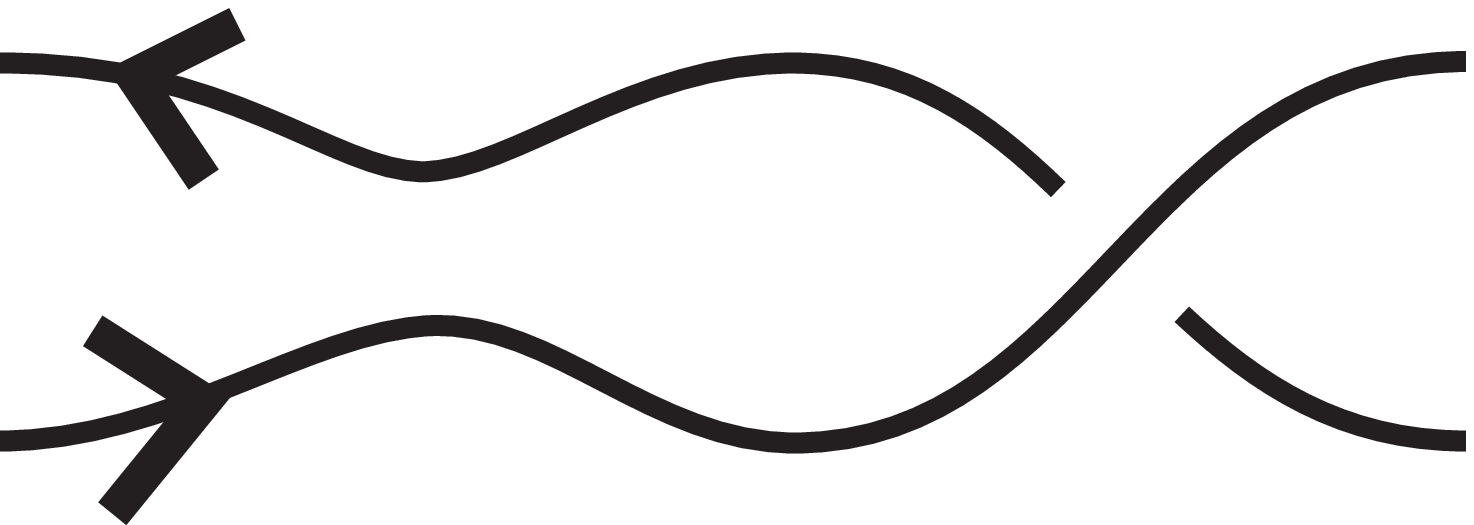}}. 
Thus we obtain a new diagram $D'$ for $K$, and we obtain the unknot by simultaneously undoing
 all the clasps.
 Put a 3-ball $B$ in the complement of $K$
 and isotope $K$ so that all the clasps are in $B$,
 and span a Seifert surface $F$ for $K$
 as in the left part of Figure \ref{fig:ooK}.
 Then $F$ decomposes into two Seifert surfaces $F_1$
 and $F_2$, where $\partial F_2$ is the unknot and
 $\partial F_1$ is the result of undoing all the clasps in $K$,
 and hence is the unknot.
\end{proof}
 
\incf{98}{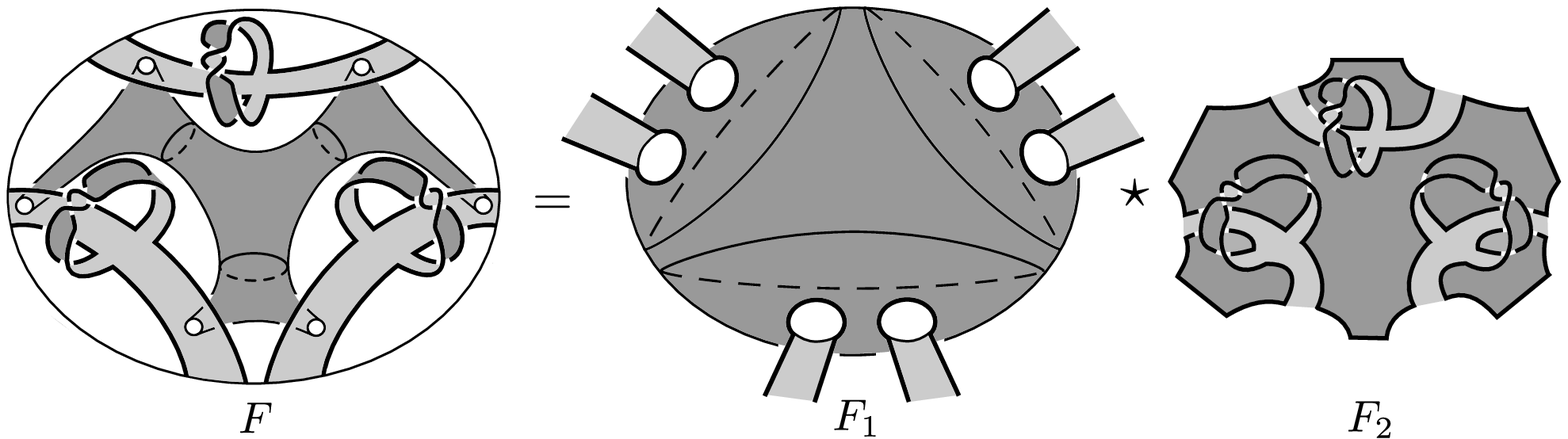}{Producing $\partial F=K$ as a Murasugi sum of unknots}{fig:ooK}

Conversely, any knot is obtained from an unknot
by simultaneously removing antiparallel clasps.
Therefore, in the proof above, we may regard
$F$ and $F_2$ as surfaces for the unknot and $F_1$ as a surface for $K$. This 
gives the following.
 
 \begin{cor}\label{cor:unknoting-by-summing-unknot}
 Any knot $K$ has a Seifert surface $F_1$ which 
 becomes a Seifert surface for the unknot by Murasugi summing $F_1$ 
with some Seifert surface for the unknot.
 \end{cor}

\begin{proof}[Proof of Theorem A1.]
By Corollary \ref{cor:unknoting-by-summing-unknot}, there is a Seifert surface $F_1$ (resp. $F_2$) for $K_1$ (resp. $K_2$) that Murasugi sums with a Seifert surface $F_1'$ (resp. $F_2'$) of the unknot $O$ to yield a Seifert surface for $O$. By Lemma \ref{lemm:sum-of-unknots}, there exist two Seifert surfaces $F_3,F_3'$ for $O$ that Murasugi sum to a Seifert surface for $K_3$. The boundary connected sum of $F_1,F_2',F_3$ (resp. $F_1',F_2,F_3'$) is a Seifert surface for $K_1$ (resp. $K_2$), and we Murasugi sum these surfaces as in Figure \ref{fig:JKL} along the shaded $n$-gon.
\end{proof}

\incf{94}{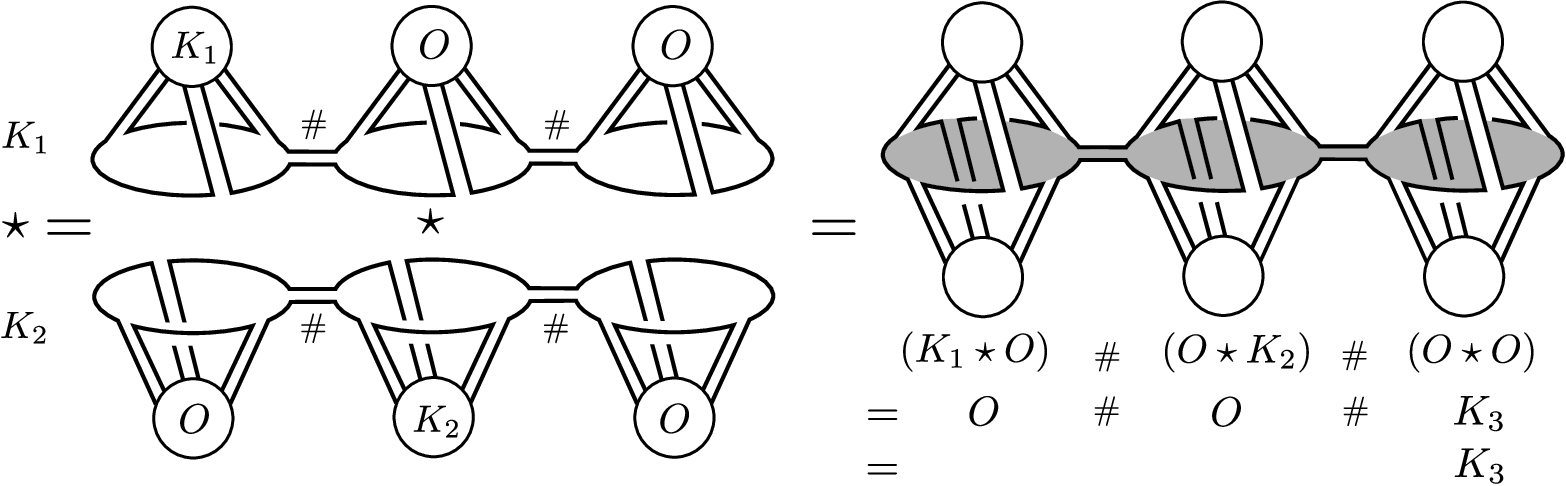}{$K_3$ as a Murasugi sum of $K_1$ and $K_2$}{fig:JKL}

\begin{proof}[Proof of Theorem A2.]
By Theorem A1, there are Seifert surfaces $F_1,F_2,F_3$ for $K_1,K_2,K_3$ 
such that $F_3$ is a Murasugi sum of $F_1$ and $F_2$ and the summing disk $\Omega$ is flat.
Apply a trivial twist at each band attached
to $\Omega$ as depicted in Figure \ref{fig:algo}.
Then we have a diagram $D'$ such that
the summing disk is spanned by a Seifert circle $C$.
Note that the canonical Seifert surface 
$F_3'$ is a Murasugi sum
of $F_1'$ and $F_2'$ along the summing disk 
$\Omega$, where 
$ \partial F_1' =\partial F_1,  \partial F_2' =\partial F_2, \partial F_3'  =\partial F_3$.
Apply Yamada's braiding algorithm \cite{yam} to $D'$, 
independently inside and outside $C$.
Then we have the desired braids. 
\end{proof}

\incf{90}{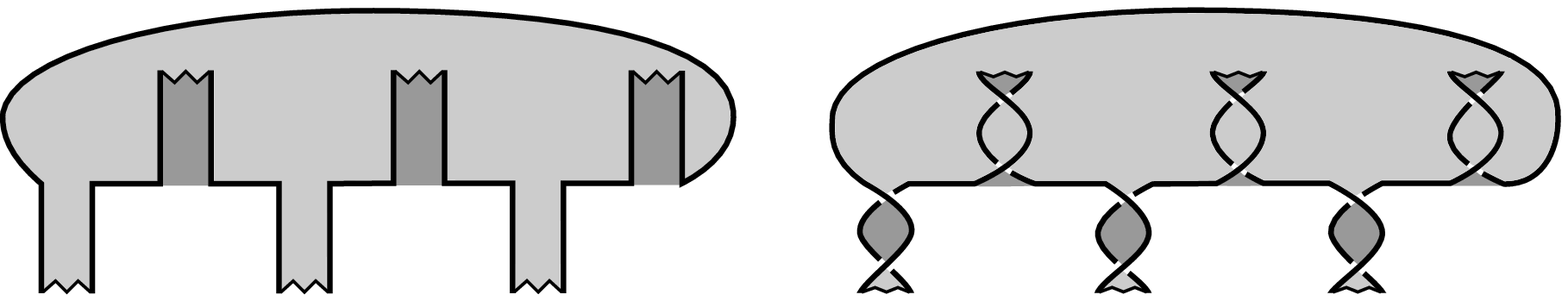}{Making a canonical surface for the braid decomposition}{fig:algo}

We use the notation $F\xrightarrow{\partial}K$ to mean that $\partial F=K$, and $F_1\stackrel{\partial}{=}F_2$ to mean that $\partial F_1=\partial F_2$.

\begin{exmp}\label{ex:knots-braid}
In Figure \ref{fig:ex1}, we illustrate $7_5$ as a Murasugi sum of two unknots, and in Figure \ref{fig:ex2}, we illustrate the unknot as a Murasugi sum of $5_2$ and $7_5$.  In these figures, we also express the Murasugi sums in terms of braid decompositions, as described in Theorem A2.  Note that we simplified some procedures in the proofs of Theorems A1 and A2.
In particular, we already have a Seifert circle corresponding to the summing disk without twisting as in Figure 6,
we eliminated some Seifert circles with two bands, and
we used ``long" bands to save us from depicting many generators.
\end{exmp}

\incf{84}{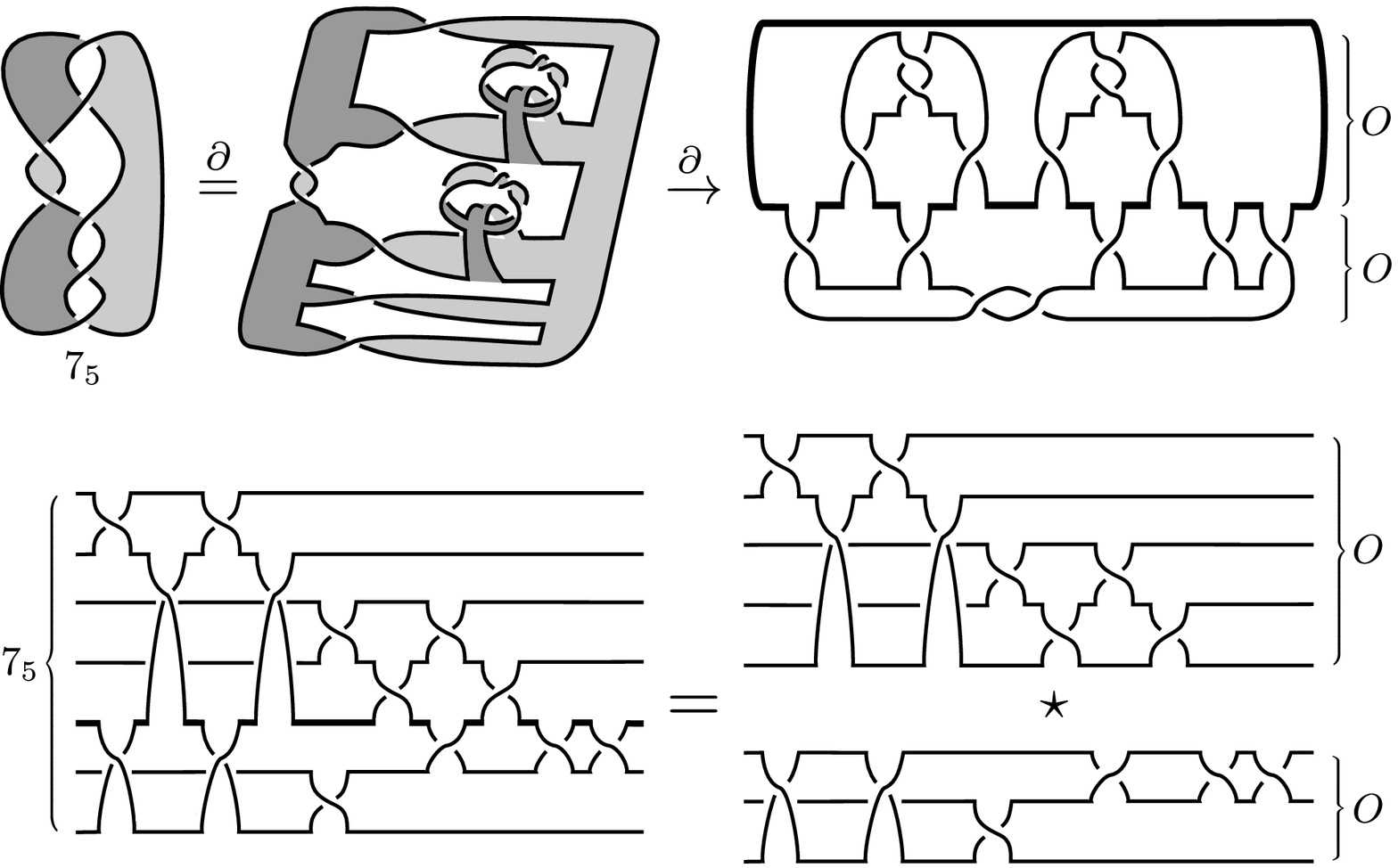}{$7_5$ as a Murasugi sum of two unknots, and its braid decomposition}{fig:ex1}

\incf{91}{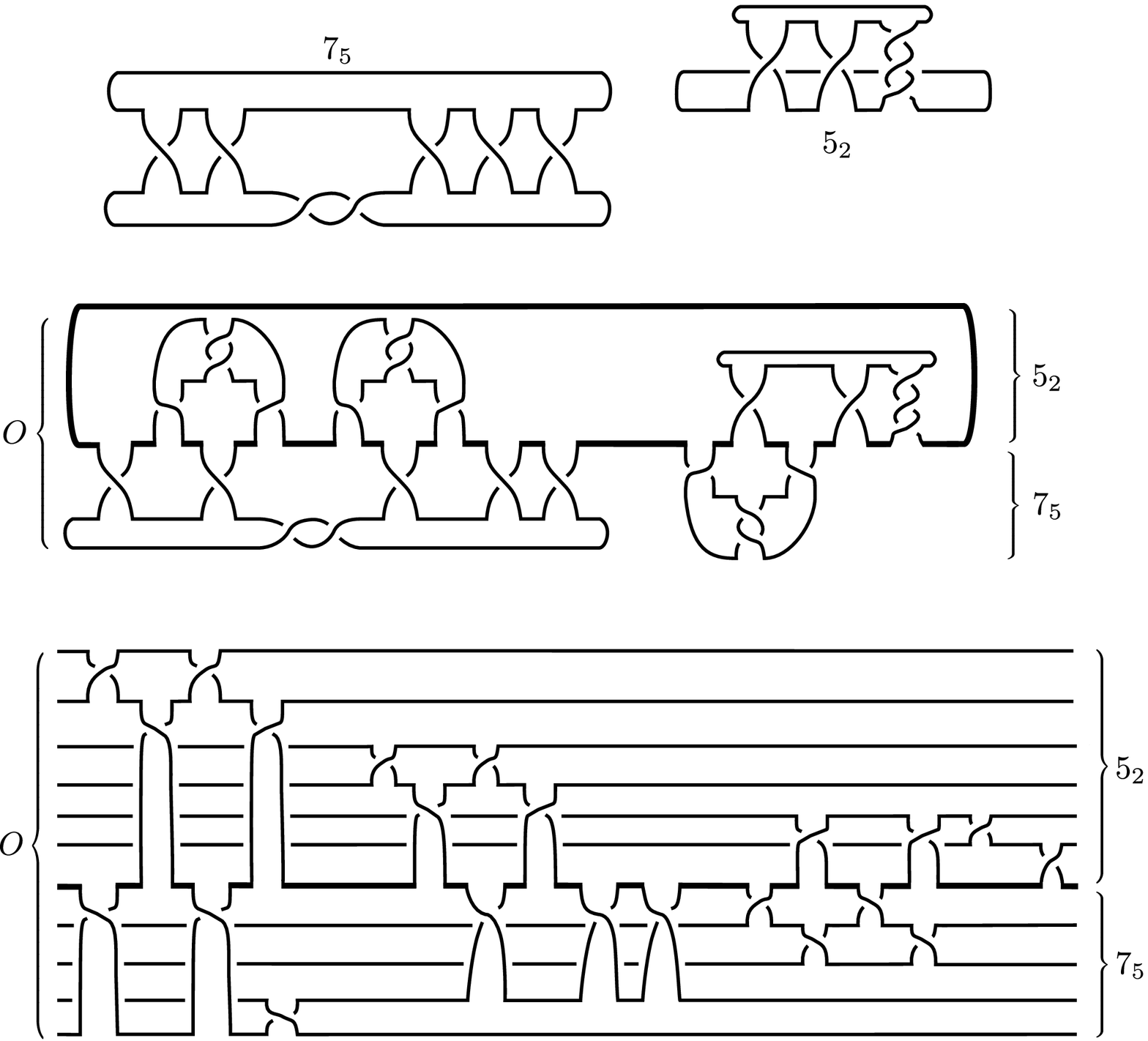}{The unknot as a Murasugi sum of $5_2$ and $7_5$, and its braid decomposition}{fig:ex2}

\section{Constraining the minimal $m$-gon}\label{sec:ngon}

In Theorem A1, it was shown that given three knots $K_1,K_2,K_3$, 
we can form $K_3$ as a Murasugi sum of $K_1$ and $K_2$ by using a sufficiently large $m$-gon. 
We begin this section with Definition \ref{dfn:m-distance} in order to give lower and upper bounds for the minimal such $m$-gon, culminating in Theorem B from the Introduction. In particular, if  we restrict the size of our $m$-gons, Theorem B obstructs forming a knot as a Murasugi sum of two given knots. Previously, the only restriction on Murasugi sums was for plumbings of fiber surfaces, such as in \cite{melmor}, where Melvin and Morton showed that the Conway polynomials of fibered knots of genus 2 take a restricted form when the fiber surface is a plumbing of Hopf bands.

\begin{defn}
\label{dfn:m-distance}
Given three oriented knots $K_1,K_2,K_3$, define the 
\emph{minimal size of Murasugi summation} for $K_3=K_1\star K_2$ to be
\begin{center}
$d_M(K_1,K_2;K_3)=\min\{m\;|\; F_3=F_1\star_m F_2, \; \partial F_i=K_i\text{ for }i=1,2,3\},$
\end{center} 
where the minimum is taken over $F_1,F_2$.
\end{defn}

Recall that $m$ is even. From the definition, $d_M(K_1,K_2;K_3)=2$ if and only if $K_3=K_1\# K_2$. 
Also, if $K_3 \neq K_1\# K_2$ can be realized as a plumbing of $K_1$ and $K_2$, then $d_M(K_1,K_2;K_3)=4$. 
Now we recall the notion of a coherent band surgery. See Figure \ref{fig:bandsurgex}.

\incf{85}{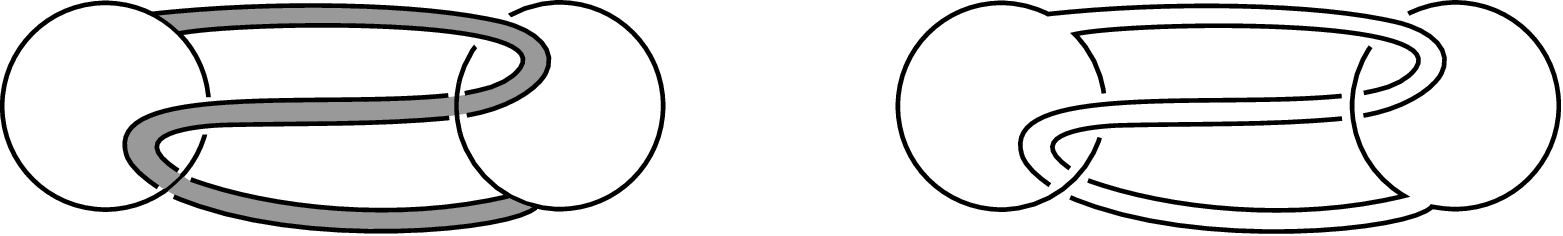}{A band surgery along an unlink yielding $3_1\# \overline{3_1}$}{fig:bandsurgex}

\begin{defn}\label{dfn:bandsurg}
Given a knot or link $L$ and an embedded band $b:I\times I\rightarrow S^3$ with $L\cap b(I\times I)=b(I\times \partial I)$, we obtain a new link $L'=(L\setminus b(I\times \partial I))\cup b(\partial I\times I)$, and we say that $L'$ is obtained from $L$ by a \emph{band surgery}. For oriented $L,L'$, a \emph{coherent} band surgery is a band surgery that respects the orientations of both links, that is, $L\setminus b(I\times \partial I)=L'\setminus b(\partial I\times I)$ as oriented spaces.
 \end{defn}

Given two oriented links $L,L'$, we denote by $d_{cb}(L,L')$ the minimal number of coherent band surgeries required to produce $L'$ from $L$. This number is known as the coherent band-Gordian distance, and in the case of knots it is equal to twice the $SH(3)$-Gordian distance \cite{kan}. In \cite{kanmor}, $d_{cb}$ is calculated for most pairs of knots up to seven crossings. For knots $K,K'$, note that $d_{cb}(K,K')$ is necessarily even because a coherent band surgery changes the number of components by one.

We have the following lower bounds for $d_M$ in terms of $d_{cb}$ and the signature $\sig$.

\begin{thm}\label{thm:lowerbounds}
For knots $K_1,K_2,K_3$, we have
\begin{center}
$ d_{cb}(K_1\# K_2,K_3)+2\leq d_M(K_1,K_2;K_3).$
\end{center}
Consequently,
\begin{center}
$ d_{cb}(K_1\sqcup K_2,K_3)+1\leq d_M(K_1,K_2;K_3),$
\end{center}
where $K_1 \sqcup K_2$ is a split link, and
\begin{center}
$|\sig(K_1)+\sig(K_2)-\sig(K_3)|+2\leq d_M(K_1,K_2;K_3).$
\end{center}
\end{thm}
\begin{proof}
Suppose $K_3$ is an $m$-Murasugi sum of $K_1$ and $K_2$, where $m=2n$ is minimal. Then there is a sequence of $2(n-1)=m-2$ coherent band surgeries between $K_3$ and $K_1\#K_2$, where each band lies within a Seifert surface for $K_3$, so that
\begin{center}
$d_{cb}(K_1\#K_2,K_3)+2\leq m.$
\end{center}
See Figure \ref{fig:obstruct}. With one more band surgery, we have $m-2+1$ coherent band surgeries between $K_3$ and the split link $K_1 \sqcup K_2$, so that 
\begin{center}
$d_{cb}(K_1\sqcup K_2,K_3)+1\leq m.$
\end{center}
If an oriented link $L'$ is obtained from $L$ by a coherent band surgery, an estimate was given by Murasugi \cite{mura2} on the difference of the signatures as $|\sig(L)-\sig(L')|\leq 1$. Since $\sig(K_1\# K_2)=\sig(K_1)+\sig(K_2)$, we obtain the third inequality. 
\end{proof}

\incf{65}{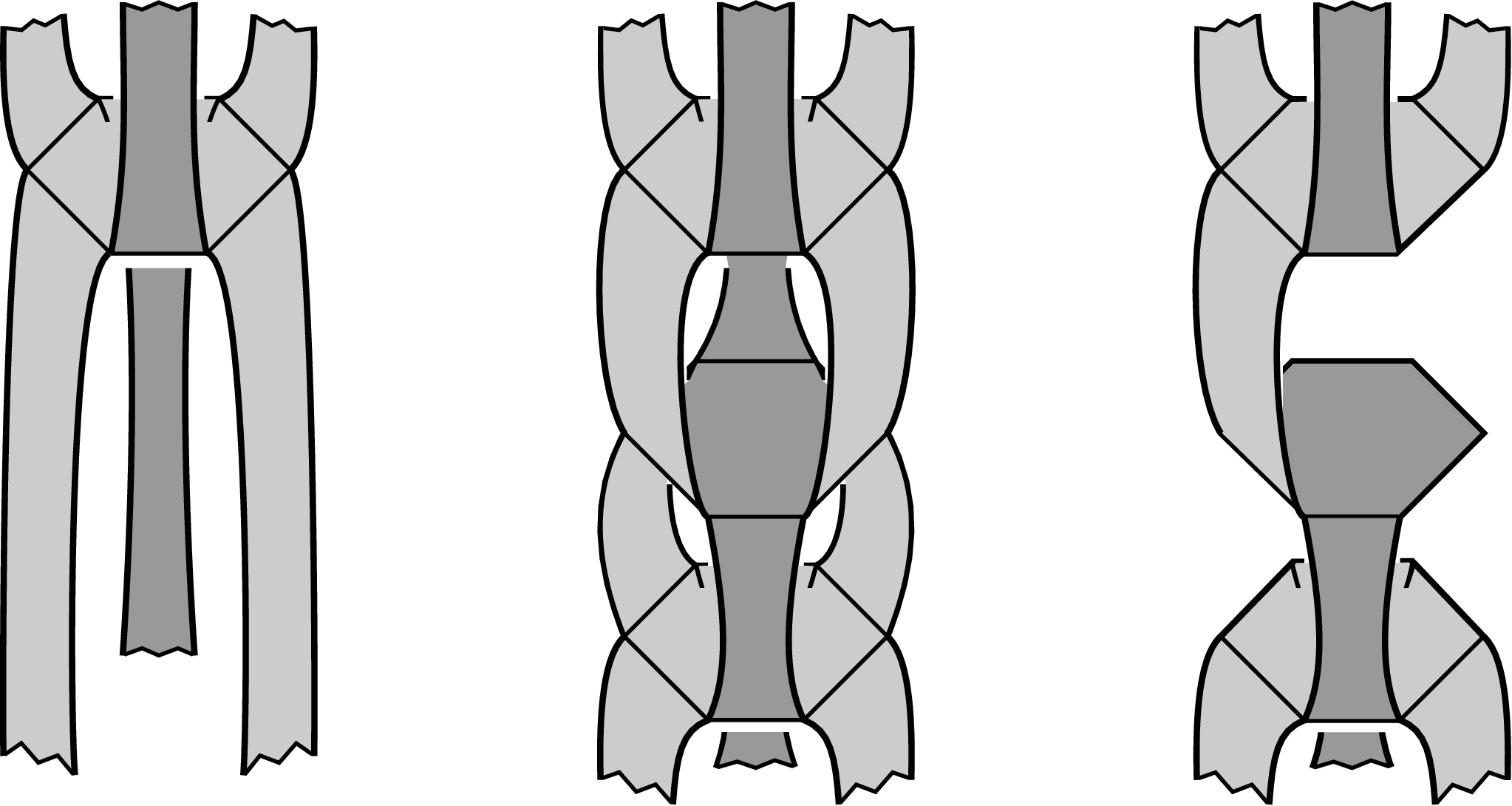}{Performing $(6-2)$ band surgeries to recover $K_1\#K_2$ from $K$}{fig:obstruct}

\begin{exmp}
By Theorem \ref{thm:lowerbounds}, the signature obstructs forming $9_1$ as a plumbing of two copies of $3_1$. However, Figure \ref{fig:9_1=3_1+3_1} depicts how $9_1$ desums into two copies of $3_1$ along the shaded 6-gon, which is obtained by first merging two 4-gons into an 8-gon which is then reduced to the 6-gon. This process is explained in the proof of Lemma \ref{lemm:mergesum}.
\end{exmp}

\incf{93}{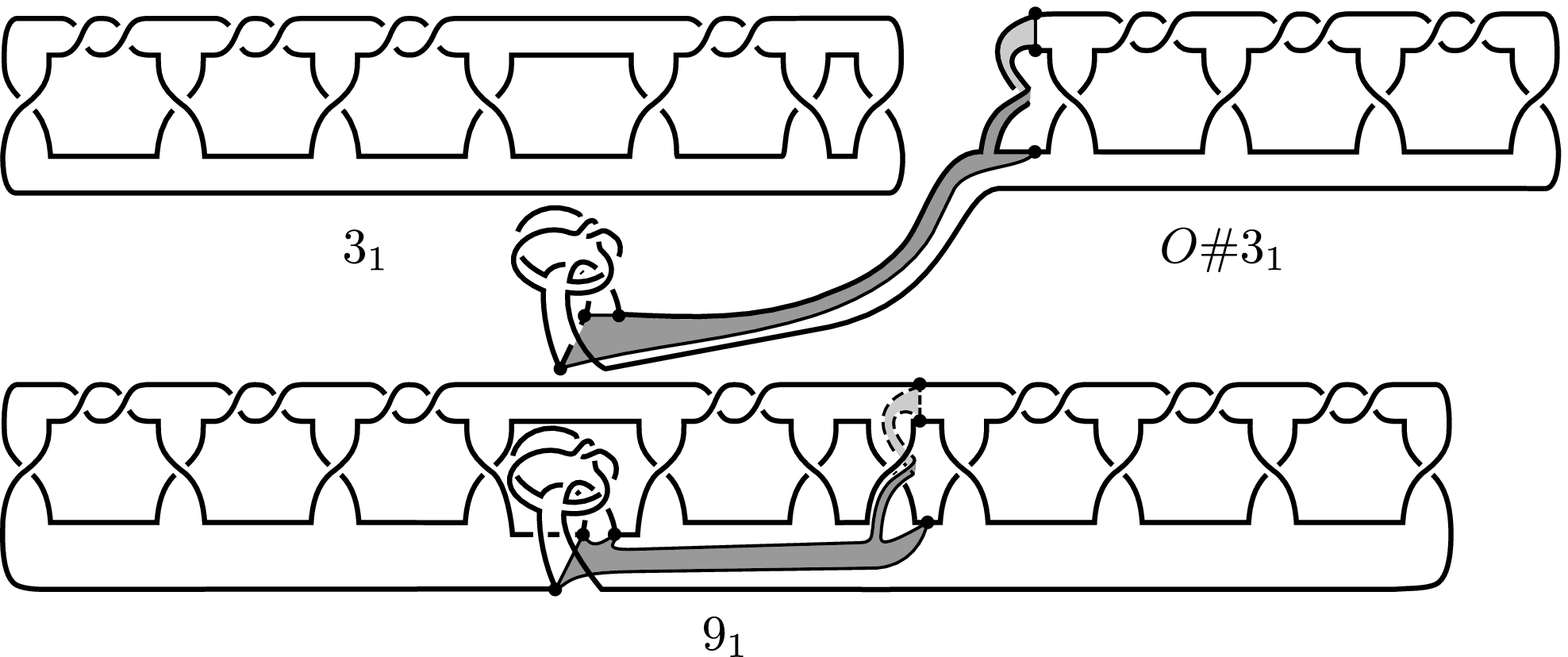}{$9_1$ as a 6-Murasugi sum of $3_1$ and $3_1$}{fig:9_1=3_1+3_1}

\begin{rmk}
For pairs of knots where $d_{cb}$ has not been determined, we may apply lower bounds for $d_{cb}$ from \cite{kan} in terms of the smooth four-ball genus $g_4$, hence in terms of $\sig$ by \cite{mura2}, and in terms of the Nakanishi index $e$. Indeed, for knots $K_1,K_2$ we have
\begin{center}
$|\sig(K_1)-\sig(K_2)| \leq 2g_4(K_1\# - K_2) \leq d_{cb}(K_1, K_2),$
\end{center}
where $-K_2$ is the reverse mirror of $K_2$, and
\begin{center}
$|e(K_1)-e(K_2)| \leq d_{cb}(K_1,K_2).$
\end{center}
\end{rmk}

We now move on to some upper bounds of $d_M$. Using Lemma \ref{lemm:band-twists-at-once} below, we can easily modify the proof of Theorem A1 to form $K_3$ as a $4(u(K_1)+u(K_2)+u(K_3))$-Murasugi sum of $K_1$ and $K_2$, where $u(K)$ is the unknotting number of $K$. This gives the upper bound
\begin{center}$d_M(K_1,K_2;K)\leq 4(u(K_1)+u(K_2)+u(K_3)).$\end{center}

In what follows, we improve this upper bound into Theorem \ref{thm:improve}. Recall that the Gordian distance $d_G$ between two knots $K,K'$ is the 
minimal number of crossing change operations required to transform a 
diagram for $K$ into a diagram for $K'$, where the minimum is taken over all diagrams for $K$.
Since $K$ can be transformed into $K'$ by a sequence of crossing changes that passes through the unknot,
we have $d_G(K,K')\leq u(K)+u(K')$. More generally, we consider the band-twist distance $d_{bt}$.

\begin{defn}\label{dfn:bt-distance}
Define the {\it band-twist distance} between two oriented links $L,L'$, 
denoted by $d_{bt}(L,L'),$ as the minimal $n$ such that there exists a 
sequence of links $L=L_0,L_1,L_2,\dots, L_n=L'$,
where $L_{i+1}$ is obtained from $L_{i}$ by an antiparallel full-twisting (a band-twist) as in Definition \ref{dfn:bandtwist}.
\end{defn}

Since any crossing change may be realized as a band-twist, we have $d_{bt}(K,K')\leq d_G(K,K').$ Just as with crossing changes, we may perform any number of band-twist operations simultaneously.

\begin{lem}\label{lemm:band-twists-at-once}
For two links $L$ and $L'$ with $d_{bt}(L,L') =n$, 
there exists a Seifert surface $F$ for $L'$
with a set $A$ of $n$ mutually disjoint properly embedded arcs such that
a Seifert surface $F$ for $L$ is obtained by applying a band-twist operation
along each arc in $A$.
\end{lem}

\begin{proof}
Let $L=L_0,L_1,L_2,\dots, L_n=L'$ be a sequence of
links related by the band-twisting operation.
After obtaining $L_1$ from $L_0$ by twisting along an
arc $\alpha_0$, instead of erasing $\alpha_0$, isotope it so that
it is disjoint from the arc $\alpha_1$ used to obtain $L_2$.
Repeating this, we obtain a set of arcs $\alpha_0, \alpha_1, \dots, \alpha_{n-1}$
attached to $L_n=L'$.
By an isotopy, we can arrange the arcs to be short and contained in a ball $B$ with
$L'\cap B$ being a trivial $2n$-string tangle.
Splice $L'$ along the arcs and push the resulting link $L''$ slightly off $B$.
Span a Seifert surface for $L''$ disjoint from $B$. 
Then we obtain the desired Seifert surface for $L'$ by attaching bands to $L''$
that pass through $B$.
\end{proof}

If we wish to perform several Murasugi sums of several surfaces with a single surface, we can often combine these sums into a single sum as indicated by the following lemma. One implication of the below is that we may perform any number of band-twist operations simultaneously via a single Murasugi sum.

 \begin{lem}\label{lemm:mergesum}
 Suppose that a Seifert surface $F'$ is
 obtained by Murasugi summing
 $S_1, S_2, \dots, S_n$
 on the same side of a connected Seifert surface $F$ along mutually disjoint 
 summing disks $\Om_1, \Om_2, \dots, \Om_n$.
 Let $\Om_i$ be an $e_i$-gon $(i= 1, 2, \dots, n)$.
 Then $F'$ is a Murasugi sum of $F$ with
 a boundary connected sum of $S_1, S_2, \dots, S_n$
 along an $m$-gon, where  $m=\sum_{i=1}^n e_i$.
 Moreover,
  if $F$ is a Seifert surface for a knot,
 then we can merge the sums into an
 $m'$-gonal sum with
  $m'=m-2(n-1)$.
 \end{lem}

\begin{proof}
We may assume that each summand $S_i$ is contained in a thin blister neighborhood of
 the summing disk $\Om_i$.
Denote the edges in  $\partial \Om_i$ as $a_{i,1},b_{i,1},a_{i,2}, b_{i,2}, \dots$ for $i=1,2,\dots,n$,
where the $a_{i,\cdot}$'s are sub-arcs of $\partial F$ and the $b_{i,\cdot}$'s are properly embedded arcs in
$F$.
Since $F$ is connected, we may find an embedded arc $\gamma$ in $F$ whose endpoints are
the midpoints of $b_{1,j}$ and $b_{2,k}$ for some $j$ and $k$ as in Figure \ref{fig:mergedisks}.
We merge the two (dark shaded) summing disks $\Omega_1$ and $\Omega_2$ into an $(e_1+e_2)$-gon 
$\Omega'$ whose boundary consists of
$(\left(\partial \Omega_1 \cup \partial \Omega_2)\setminus (b_{1,j}\cup b_{2,k})\right) \cup (\gamma_1 \cup \gamma_2)$,
where $\gamma_1$ and $\gamma_2$ are properly embedded in $F$ and
$b_{1,j} \cup \gamma_1 \cup b_{2,k} \cup \gamma_2$ is a rectangle $R=\gamma \times I$ such that 
$R \cap {\rm int}\! \bigcup_{i=1}^n  \Omega_i$=$\emptyset$.
We see that the boundary connected sum of $S_1$ and $S_2$ is contained in a thin blister neighborhood of the new summing disk $\Omega'$.
By repeating this merging operation, we eventually combine all the summing disks into one and
obtain the desired Murasugi sum.

For the last part of the assertion, the assumption that $\partial F$ being connected ensures the
existence of two arcs $a_{1,p}, a_{2,q}$ for some $p, q$ such that one segment of $\partial F$ between
$a_{1,p}$ and $a_{2,q}$ does not pass through other summing disks.
Then we can apply the previously mentioned merging of $\Omega_1$ and $\Omega_2$ so that
$\gamma_1$ or $\gamma_2$, 
say $\gamma_1$, can be isotoped to a sub-arc of 
$\partial F$ in $F \setminus \bigcup_{i=1}^n \Om_i$.
Then the three consecutive edges in $\Omega'$, 
say, $a_{1,k}, \gamma_1, a_{2,j}$ for some $k,j,$
can be regarded as one edge by merging the bands
of $S_1$ and $S_2$ attached to $a_{1,k}$ and $a_{2,j}$ as in Figure \ref{fig:reducedisks}.
Then the new summing disk $\Omega''$ 
is a $p$-gon, where $p=e_1+e_2 -2$.
\end{proof}

\incf{70}{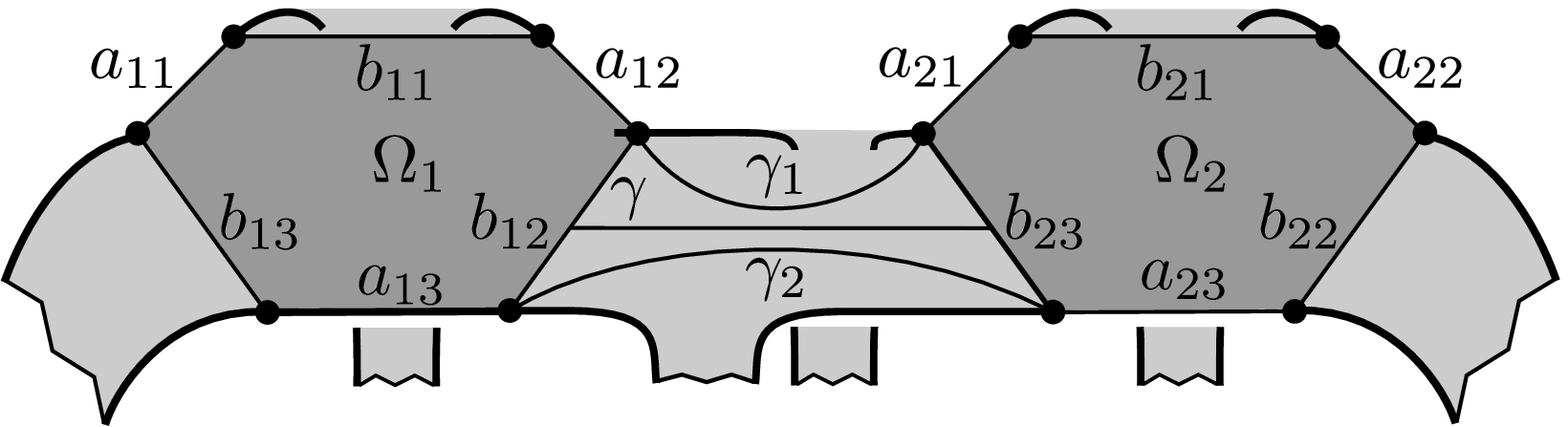}{Merging disks $\Om_1,\Om_2$ in $F$ along $\g$}{fig:mergedisks}

\incf{99}{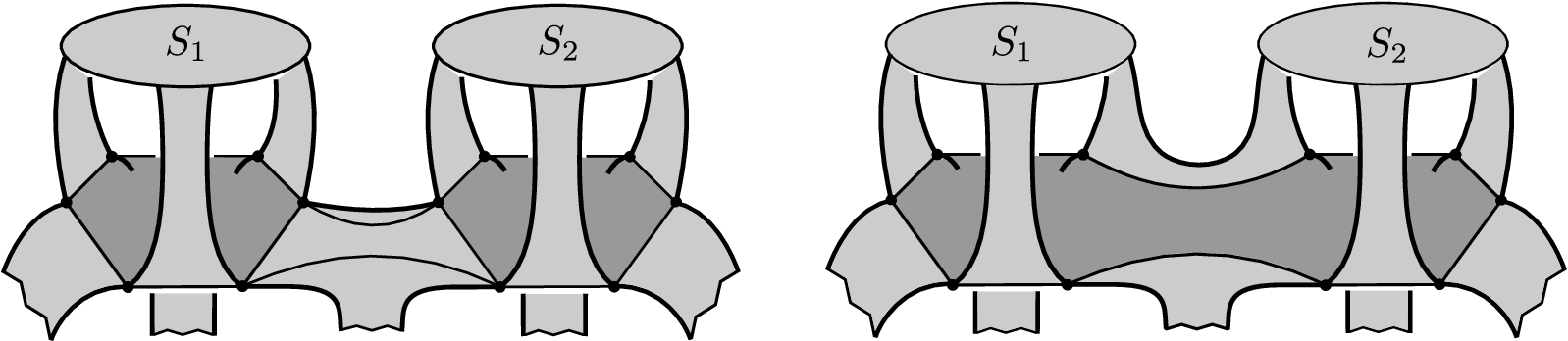}{Merging two 6-gons into a ($6+6-2$)-gon}{fig:reducedisks}

Combining Lemma \ref{lemm:band-twists-at-once} and Lemma \ref{lemm:mergesum}, we arrive at the following improvement on our upper bound.

\begin{thm}\label{thm:improve}
Let $K_1,K_2,K_3$ be knots. Then 
\begin{center}$d_{M}(K_1,K_2;K_3)\leq 2(d_{bt}(K_1,K_3)+d_{bt}(K_2,O)+1),$\end{center}
where the roles of $K_1$ and $K_2$ may be switched to improve the upper bound. 
\end{thm}

\begin{proof}
Suppose that $d_{bt}(K_1,K_3)=p$ and $d_{bt}(K_2,O)=q$. As guaranteed by Lemma \ref{lemm:band-twists-at-once}, there is a Seifert surface $F_1$ for $K_1$ (resp. $F_2$ for $K_2$) with a collection $\A$ of arcs $\alpha_1,\ldots, \alpha_p$ (resp. $\B$ of arcs $\beta_1,\ldots,\beta_q$) such that performing band-twist operations along the arcs yields a Seifert surface $F_3$ for $K_3$ (resp. $F_0$ for $O$). Prepare Seifert surfaces $A_1,\ldots, A_p$ and $B_1,\ldots, B_q$ such that plumbing $A_j$ along $\alpha_j$ (resp. $B_k$ along $\beta_k$) results in applying the band-twist operation along $\alpha_j$ (resp. $\beta_k$). More precisely, each of the surfaces is a plumbing of a trivial annulus and an unknotted annulus with various numbers of full-twists (recall Figure 2).

Construct $F_1'$ from $F_1$ such that $\partial F_1'=K_3$ by plumbing $A_1,\ldots, A_p$ along $\A$ on the positive side of $F_1$. Also, construct $F_0'$ from $F_0$ such that $\partial F_0'=O$ by plumbing $B_1,\ldots, B_q$ along $\B$ on the negative side of $F_0$. 

We merge the plumbed surfaces $A_1,\ldots A_p$ so that $F_1'$ is a Murasugi sum of $F_1$ and $A$, where $A$ is a boundary connected sum of $A_1,\ldots, A_p.$ Note that $\partial A$ is an unknot $O_1$. By Lemma \ref{lemm:mergesum}, we may regard the Murasugi sum as along a $(4p-2(p-1))$-gon and hence a $(2p+2)$-gon. Similarly, $F_0'$ is regarded as a $(2q+2)$-Murasugi sum of $F_0$ and a Seifert surface $B$ for an unknot $O_2$. Denote by $F$ the 2-Murasugi sum (i.e., a boundary connected sum) of $F_1'$ and $F_0'$, where $F_0'$ is summed on the positive side of $F_1'$. 

Then $F$ is a $(2p+2+2q+2)$-Murasugi sum of two summands, where one summand is the boundary connected sum of $F_1$ and $B$, and the other summand is the boundary connected sum of $F_0$ and $A$. The boundary of the first (resp. second) summand is $K_1\#O_2$ (resp. $K_2\#O_1$). By Lemma \ref{lemm:mergesum}, we can reduce the summing $(2p+2q+4)$-gon into a $(2p+2q+2)$-gon. Therefore, we have expressed $K_3$ as a $2(p+q+1)$-Murasugi sum of $K_1$ and $K_2$.
\end{proof}

By combining Theorems \ref{thm:lowerbounds} and \ref{thm:improve}, we arrive at Theorem B. As an application of Theorem B, we determine $d_M(3_1,3_1;K)$ for knots up to five crossings. Theorem B shows that
\begin{center}$d_M(3_1,3_1;3_1)=4, \qquad d_M(3_1,3_1;O)=d_M(3_1,3_1;4_1)=6,$\end{center}
while it only gives the bounds 
\begin{center}$4\leq d_M(3_1,3_1;5_1),\; d_M(3_1,3_1;5_2) \leq 6.$\end{center}
We can show that $d_M(3_1,3_1;5_1)=d_M(3_1,3_1;5_2)=4$ in the following way.

For $a_1,a_2,\ldots, a_n\in 2\bbz$, denote by $S[a_1,a_2,\ldots,a_n]$ a linear plumbing of $n$ unknotted annuli, where the $i^{th}$ annulus has $a_i$ half-twists. Note that all 2-bridge links have such a linear plumbing as a Seifert surface \cite{hatthu}, which is of minimal genus if and only if $a_1a_2\cdots a_n\neq 0$. Using the notation of Example \ref{ex:knots-braid}, we have the following:
\begin{enumerate}
\item $S[a_1,a_2,\ldots,a_n]\star_4 S[b_1,b_2,\ldots,b_m] =S[a_1,\ldots,a_n,b_1,\ldots,b_m]$,
\item $S[a_1,a_2,\ldots,a_n]\stackrel{\partial}{=}S[a_1,a_2,\ldots,a_n,a_{n+1},0]$ for any $a_{n+1}$,
\item $S[a_1,a_2,\ldots,a_i,0,a_{i+2},\ldots, a_n]\stackrel{\partial}{=}S[a_1,a_2,\ldots,a_i+a_{i+2},\ldots, a_n]$.
\end{enumerate}
Using the above, we have:

\begin{exmp}\label{thom3141}

\begin{enumerate}
\item $3_1\star_4 3_1 \xleftarrow{\partial} S[2,2]\star_4 S[2,2]=S[2,2,2,2]\xrightarrow{\partial}5_1$
\item $3_1\star_4 3_1  \xleftarrow{\partial} S[2,2,-2,0]\star_4 S[2,2]\stackrel{\partial}{=}S[2,2,0,2]\stackrel{\partial}{=}S[2,4]\xrightarrow {\partial}5_2.$
\end{enumerate}
\end{exmp}

Using this notation, we summarize what Thompson showed in \cite{thomp} as follows:

\begin{enumerate}
\item $O\star_4 O\xleftarrow{\partial}S[2,0]\star_4 S[0,2]=S[2,0,0,2]\stackrel{\partial}{=}S[2+0,2]=S[2,2]\xrightarrow{\partial}3_1$
\item $4_1 \star_4 4_1 \xleftarrow{\partial} S[2,-2,2,0]\star_4 S[-2,2]\stackrel{\partial}{=} S[2,-2,0,2]\stackrel{\partial}{=}S[2,0]\xrightarrow{\partial}O.$ 
\end{enumerate}

\section{Concluding remarks and further discussion}

When Thompson \cite{thomp} negated the possibility of generalizing nice properties
preserved under Murasugi sums and decompositions, she also initiated the study of
constructing knots by Murasugi sums with her examples of $3_1 = O \star_4 O$ and $O=4_1 \star_4 4_1$.
Around that time, it was shown in \cite{haywad} that any link is constructed by successively plumbing
trivial (i.e., unknotted and untwisted) annuli, where the gluing square could be quite complicated.
Then in \cite{fhk}, it was shown that any link has a Seifert surface which is obtained from a disk $D$
by successively plumbing trivial annuli to $D$, where all of the gluing squares are in $D$,
thus giving a new standard form of links called a {\it flat basket plumbing presentation}.

In this paper, we have shown that given a knot $K$ and any two knots $K_1$, $K_2$,
we can form $K$ as a Murasugi sum of $K_1$ and $K_2$, and that this situation can be illustrated in
a closed braid form.
So anything can happen when we use Murasugi sums of general complexity, but as we have also shown,
Murasugi sums are more well-behaved when the size of Murasugi sums is restricted.
For further study, we may also put other restrictions on the Murasugi sums, for example by
restricting the genera of the Seifert surfaces involved in the Murasugi sums.
So we can ask the following decomposition question about the result of the Murasugi sum:
\begin{question}
For a nontrivial knot $K$ with genus $g(K)$, what is the minimal genus
of a Seifert surface $F$ for $K$ which is a Murasugi sum of two unknots?
\end{question}
In this situation, if $K$ has an alternating diagram which can be unknotted with $u(K)$ crossing changes, then Seifert's algorithm yields a minimal genus Seifert surface, so our constructions give that $1\leq g(F)-g(K)\leq u(K)$. Note that this inequality is also true for knots with $u(K)=1$, because by \cite{koba}, there exists a minimal genus Seifert surface for $K$ on which the
unknotting crossing change can be done by twisting a band. Similarly, we can ask the following composition question about the summands of the Murasugi sum:
\begin{question}
For two nontrivial knots $K_1,K_2$, what are the minimal genera of Seifert surfaces $F_1,F_2$ of $K_1, K_2$ that Murasugi sum to the unknot?
\end{question}
Answering these questions in general seems to be more involved than what we treat in this paper, where we typically form some Seifert surfaces, manipulate them, then dissolve them to obtain their knot boundaries. In particular, one would need to be more explicit with how the original Seifert surfaces are formed. Studying these questions via band surgery and band twisting (or perhaps other moves) should yield insight into a collection of subgraphs of $M\kern-1pt SG(K)$ similar to $M\kern-1pt SG(K,n)$.

\section{Appendix}\label{sec:append}

This section is not directly concerned with the main theorems, 
but we give proofs for Theorems C1 and C2 here to record some other ways of constructing knots and links by Murasugi sums.

\begin{proof}[Proof of Theorem C1.]
Let $\widetilde{F}$ be a canonical surface on 
a closed braid for $L$.
For each set of bands connecting the same pair of
disks, select one to save, and apply a modification as in 
Figure \ref{fig:pHopf} for the other bands. This produces a Seifert surface $F$.
Then at each site of modification, we can deplumb
a positive Hopf band. Note that
the Hopf bands are plumbed on the same side, and
hence we can merge them into a boundary connected sum of 
positive Hopf bands.
Deplumbing all Hopf bands yields a surface $F'$
which is composed of 
\begin{enumerate}
\item the Seifert disks of $F$,
\item one band between each pair of adjacent disks, and 
\item one tube at each site of the bands which are cut by the deplumbing.
\end{enumerate}
Slide the tubes to the top of the disks and we have
the desired once-punctured Heegaard surface.
\end{proof}

\incf{90}{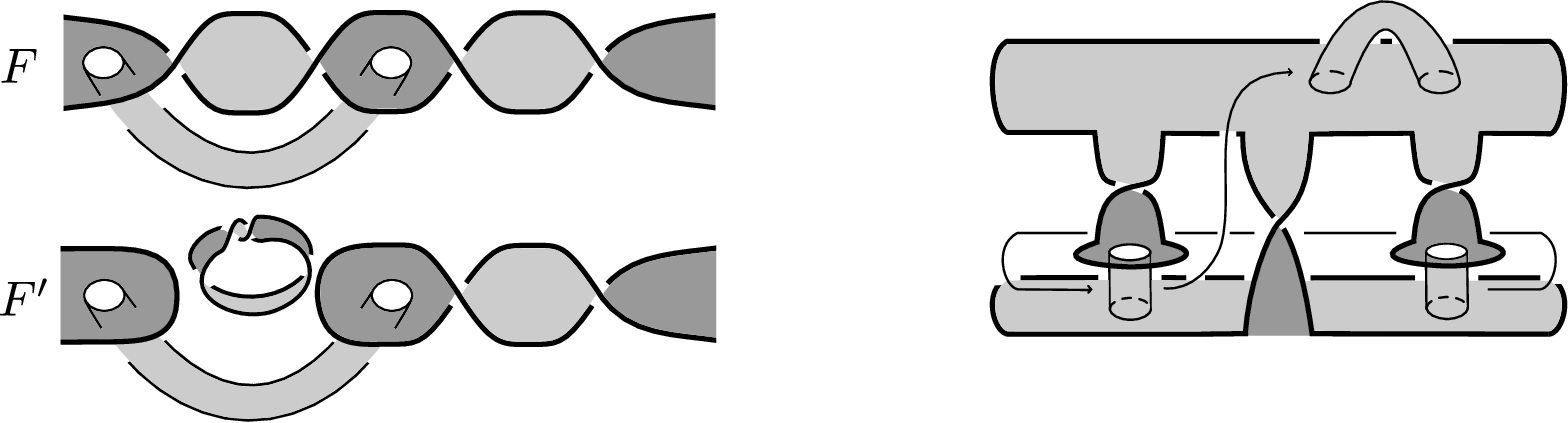}{Deplumbing a positive Hopf band, and then sliding tubes to the top disk}{fig:pHopf}

\begin{proof}[Proof of Theorem C2.]
Take a diagram $D$ of $L$, apply a Reidemeister II move at each crossing of $D$
to triple the crossing number
as in Figure \ref{fig:depflata}, and span
a canonical surface $F$.
Then we can deplumb a flat annulus at the site 
of each original crossing.  
We see that all annuli are plumbed to $F$ on the
positive side, which has lighter shading. Merge all plumbings so that
$F$ is a Murasugi sum of a Seifert surface $F'$ and a surface $P_1$, where
$P_1$ is a boundary connected sum of flat annuli, and hence a planar surface.
Now, $F'$ is obtained from Seifert disks of the original
diagram $D$ by applying a tube on the negative side of the site of each crossing. 
Therefore, $F'$ is the boundary of the
neighborhood (i.e., a standard handlebody) 
of the Seifert graph $G$ of $D$ with
punctures. In other words,
place a punctured sphere
at each vertex of $G$ and apply a tube along each 
edge of $G$.
Hence $F'$ is a standard Heegaard surface with punctures, which is isotopic to
a linear plumbing of untwisted, unknotted annuli with punctures.
We can regard $F'$ as a planar surface $P$
with
a planar surface $P_2$
Murasugi-summed on the negative side.
\end{proof}

\incf{90}{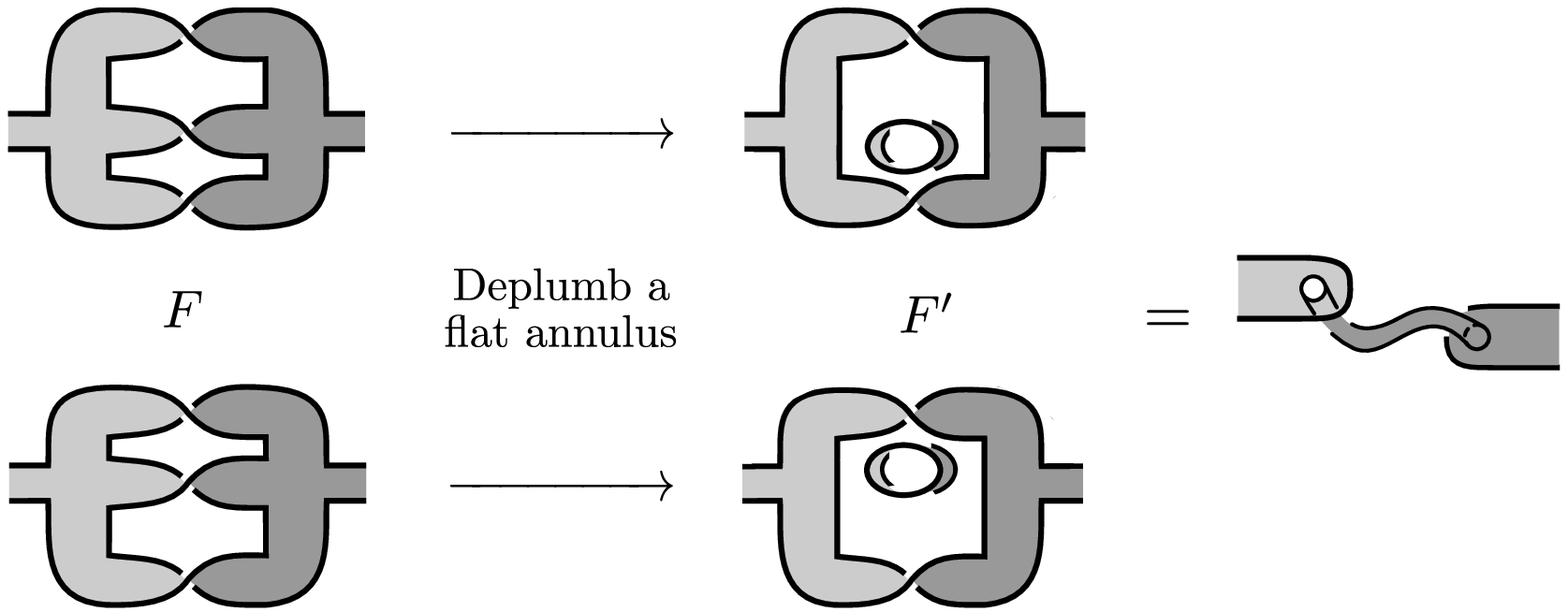}{Applying a Reidemeister II move, and then deplumbing a flat annulus}{fig:depflata}

\bibliographystyle{amsxport}  
\bibliography{./newbib2.bib}

\end{document}